\documentclass[a4paper,10pt]{article}
\usepackage{amssymb,amsmath,amsthm}
\usepackage{fullpage}
\usepackage{hyperref}
\usepackage{eucal}
\usepackage{graphicx}
\usepackage{mathrsfs}
\newcommand{\ie}{\emph{i.e.}}
\newcommand{\eg}{\emph{e.g.}}
\newcommand{\cf}{\emph{cf.}}
\newcommand{\Real}{\mathbb{R}}

\newcommand{\Sphere}{\mathbb{S}}

\newcommand{\dist}{\mathop{\mathrm{dist}}\nolimits}

\newcommand{\divergence}{\mathop{\mathrm{div}}\nolimits}
\newcommand{\diag}{\mathrm{diag}}

\newcommand{\der}{\mathrm{d}}
\newtheorem{Theorem}{Theorem}

\newtheorem{Proposition}{Proposition}

\theoremstyle{definition}
\newtheorem{Remark}{Remark}

\usepackage[normalem]{ulem}
\usepackage{color}
\definecolor{DarkBlue}{rgb}{0,0.1,0.7}
\definecolor{DarkGreen}{rgb}{0,0.5,0.1} 

\newcommand\soutD{\bgroup\markoverwith
{\textcolor{DarkGreen}{\rule[.01ex]{2pt}{1pt}}}\ULon}
\newcommand{\Hm}[1]{\leavevmode{\marginpar{\tiny%
$\hbox to 0mm{\hspace*{-0.5mm}$\leftarrow$\hss}%
\vcenter{\vrule depth 0.1mm height 0.1mm width \the\marginparwidth}%
\hbox to
0mm{\hss$\rightarrow$\hspace*{-0.5mm}}$\\\relax\raggedright #1}}}

\definecolor{grey}{rgb}{.7,.7,.7}  %

\begin{document}
%
\title{\textbf{\LARGE
The Cheeger constant of curved tubes
}}
\author{David Krej\v{c}i\v{r}{\'\i}k,\!$^a$ \
Gian Paolo Leonardi\,$^b$ \ and \ Petr Vlachopulos\,$^a$}
\date{\small 
\begin{quote}
\emph{
\begin{itemize}
\item[$a)$] 
Department of Mathematics, Faculty of Nuclear Sciences and 
Physical Engineering, Czech Technical University in Prague, 
Trojanova 13, 12000 Prague 2, Czechia;
david.krejcirik@fjfi.cvut.cz%
\item[$b)$]
Department of Physics, Informatics and Mathematics,
University of Modena and Reggio Emilia,
Via Campi, 213/b, 41100 Modena, Italy;
gianpaolo.leonardi@unimore.it
\end{itemize}
}
\end{quote}
29 November 2018}
\maketitle
\begin{abstract}
\noindent
We compute the Cheeger constant of spherical shells and tubular neighbourhoods 
of complete curves in an arbitrary dimensional Euclidean space.
%
%
\end{abstract}
%
 
\section{Introduction}
%
Given an open connected set $\Omega \subset \Real^d$ with $d \geq 1$,
we define the \emph{Cheeger constant} of~$\Omega$ to be the number
\begin{equation}\label{Cheeger}
  h(\Omega) := \inf_{S \subset \Omega} \frac{|\partial S|}{|S|}
  \,,
\end{equation}
where the infimum is taken over all non-empty bounded sets $S \subset \Omega$.
We denote by $|S|$ the volume (\ie\ $d$-dimensional Lebesgue measure) of~$S$. We also use the notation $|\partial S|$ for the perimeter of $S$. The perimeter functional extends in the appropriate, variational way the $(d-1)$-dimensional Hausdorff measure of the boundary of $S$ --  assuming $\partial S$ smooth, or Lipschitz -- to the whole class of Borel sets; for more details about its general definition and properties, we refer the interested reader to \cite{Maggi_book}.
Any minimiser of~\eqref{Cheeger}, 
if it exists (\eg\ if~$\Omega$ is bounded),
is called a \emph{Cheeger set} of~$\Omega$
and is denoted by $\mathcal{C}_\Omega$.	
We refer to the survey papers~\cite{Leonardi_2015,Parini}
for motivations to consider the problem~\eqref{Cheeger},
properties, history and many references.
In particular, the close connection between the Cheeger problem and the eigenvalue problem for the $p$-Laplacian when $p$ goes to $1$ is studied in \cite{Kawohl-Fridman_2003}, while in the recent contribution \cite{Kawohl-Schuricht_2015} one can find insights on the eigenvalue problem for the $1$-Laplacian.

There are very few known domains~$\Omega$ for which 
the Cheeger constant~$h(\Omega)$ can be computed explicitly.
In all dimensions, these are essentially just balls
$B_a := \{x\in\Real^d:|x|<a\}$ with $a>0$, 
for which one has		
\begin{equation}\label{ball}
  h(B_a) = \frac{d}{a} 
\end{equation}
and $\mathcal{C}_{B_a}=B_a$.
Apart from the trivial one-dimensional situation,
much more is known in the planar case.
It is easy to compute the Cheeger set of rectangles
and in fact there exists a constructive algorithm for 
finding the Cheeger constant in any convex polygon	
\cite{Kawohl-Lachand-Robert_2006}.

It was pointed out in~\cite{KP}
that there exists another large class of planar domains
for which the Cheeger constant can be computed explicitly,
namely \emph{curved strips}.
Given a smooth closed planar curve~$\Gamma$ 
and a positive number~$a$, we define the strip of radius~$a$
to be the tubular neighbourhood  
$\Omega_a := \{x \in \Real^2 : \dist(x,\Gamma) < a\}$.
If~$a$ is so small that~$\Omega_a$ ``does not overlap itself'',
the main result of~\cite{KP} says that
\begin{equation}\label{Cheeger.strip}
  h(\Omega_a) = \frac{1}{a} 
  \,,
\end{equation}
that $\mathcal{C}_{\Omega_a}=\Omega_a$
and that the former remains to hold 
for strips about unbounded complete curves, too. 
Moreover, precise estimates for the Cheeger constant
and explicit characterisation of the Cheeger set
for general strips were given in~\cite{KP} 
and further improved in~\cite{Leonardi-Pratelli_2016}. 
See also the recent papers \cite{Leonardi-Neumayer-Saracco_2017} (which contains a characterisation of the maximal Cheeger set within a Jordan domain of the Euclidean plane) and \cite{Leonardi-Saracco_2018} (where, motivated by capillarity-related issues, two nontrivial examples of minimal Cheeger sets in the plane are constructed).

It is frustrating that already a three-dimensional cube
does not admit an explicitly known Cheeger constant 
and there is no explicit analytical description
of its Cheeger set (see \cite[Open Problem~1]{Kawohl_2016}).
In general, much less is known about the Cheeger problem~\eqref{Cheeger}
in higher dimensions.

In this paper we prove that curved tubes and spherical shells are the unique Cheeger sets of themselves (in other words, they are minimal with respect to inclusion, hence we shall also refer to them as \textit{minimal Cheeger sets} from now on). The fact that a $d$-dimensional tubular neighbourhood of a closed curve is a minimal Cheeger set corresponds to a higher-dimensional version of the main result of~\cite{KP}. 
Given a closed smooth curve~$\Gamma$ in~$\Real^d$ with $d \geq 2$,
we introduce a \emph{curved tube} by 
\begin{equation}\label{tube}
  \Omega_a := \{x \in \Real^d : \ \dist(x,\Gamma) < a\}
  \,.
\end{equation}
We say that~$\Omega_a$ \emph{does not overlap itself}
if the map 
$
  \Gamma\times(0,a) 
  \ni (q,t) \mapsto q + t \, N(q)
$
induces a smooth diffeomorphism 
for any smooth normal vector field~$N$ along~$\Gamma$.
Since~$\Gamma$ is compact, this condition holds 
for all sufficiently small~$a$.
Then our first result reads as follows. 

\begin{Theorem}\label{Thm.main}
Given a closed smooth curve~$\Gamma$ in~$\Real^d$ with $d \geq 2$
and a positive number~$a$, let~$\Omega_a$ be defined by~\eqref{tube}.
If~$a$ is so small that~$\Omega_a$ does not overlap itself,
then 
\begin{equation}\label{Cheeger.tube}
  h(\Omega_a) = \frac{d-1}{a} 
\end{equation}
and $\mathcal{C}_{\Omega_a}=\Omega_a$.
\end{Theorem}

Comparing~\eqref{Cheeger.tube} with~\eqref{ball},
we see that the Cheeger constant in a $d$-dimensional 
curved tube of radius~$a$ equals the Cheeger constant
of the $(d-1)$-dimensional ball of the same radius.
It is worth noticing that
(contrary to the first eigenvalue of the Dirichlet Laplacian in~$\Omega_a$,
which highly depends on the geometry of~$\Gamma$) 
the shape of the underlying manifold~$\Gamma$
has no influence on the value of the Cheeger constant
(which can be interpreted as the first eigenvalue of
the non-linear $1$-Laplacian).

We now state our second result about spherical shells. 
\begin{Theorem}\label{Thm.shell}
Given two positive radii $r<R$, the spherical shell 
$A_{r,R} := \{x\in \Real^d:\ r<|x|<R\}$ is a minimal Cheeger set, and 
\begin{equation}\label{Cheeger.shell}
h(A_{r,R}) = d \, \frac{R^{d-1}+r^{d-1}}{R^d-r^d}\,.
\end{equation}
\end{Theorem}

This theorem extends the result
of \cite[Sec.~11, Ex.~4]{Bellettini-Caselles-Novaga_2002}
for annuli to higher dimensions.
Notice that our Theorems~\ref{Thm.main} and~\ref{Thm.shell} coincide
in the special case when $\Gamma$ is a circle and $d=2$.
It is a challenging open problem to determine the Cheeger constant
of tubular neighbourhoods of general submanifolds of~$\Real^d$.

\section{The Proofs}
%
Although~\eqref{Cheeger.tube} appears 
to be a natural generalisation of~\eqref{Cheeger.strip}
to higher dimensions, the strategy of proof that we adopt here for both Theorems \ref{Thm.main} and \ref{Thm.shell} substantially
differs from the one presented in~\cite{KP}.  
The latter was essentially based on the fact that 
taking the convex hull of a planar domain
simultaneously enlarges the area and reduces the perimeter.
Since this property is unavailable in higher dimensions,
here we replace the reasoning by an argument based on ``test vector fields''
(\cf~\cite[Prop.~1]{Grieser_2006}).
\begin{Proposition}\label{Prop.Grieser}
Let $\Omega \subset \Real^d$ be an open connected set.
If there exist a smooth vector field $V:\Omega\to\Real^d$ 
satisfying the pointwise inequalities 
\begin{equation}\label{test.field}
  |V| \leq 1
  \qquad \mbox{and} \qquad 
  \divergence V \geq c
\end{equation}
in~$\Omega$ with some $c \in \Real$, 
then $h(\Omega) \geq c$.
\end{Proposition}

Clearly, the proposition can be used to obtain lower bounds
to the Cheeger constant.
Upper bounds can be obtained more straightforwardly 
by using suitable ``test domains'' in~\eqref{Cheeger}. 

\subsection{Proof of Theorem \ref{Thm.main}: the upper bound}\label{Sec.upper}
To show that $h(\Omega_a)$ can be bounded from above
by the right-hand side of~\eqref{Cheeger.tube},
we use the whole tube~$\Omega_a$ as a test domain in~\eqref{Cheeger}.
To this aim, we first need to establish some basic facts
about the geometry of curved tubes.

Since we assume that~$\Omega_a$ does not overlap itself,
necessarily~$\Gamma$ is an embedded submanifold of~$\Real^d$. 
Let us parameterise $\Gamma \subset \Real^d$ locally by 
a smooth map $\gamma:I\to\Real^d$, where $I \subset \Real$
is an open interval. 
In particular, $\gamma(s) \in \Gamma$  for all $s \in I$.
Without loss of generality, we can assume that~$\gamma$ is unit-speed,
\ie\ $|\dot\gamma(s)|=1$ for all $s \in I$.
We denote by $e_1:=\dot\gamma$ and $\kappa := |\ddot\gamma|$, 
the tangent vector field and curvature of~$\Gamma$, respectively. 

Note that we allow curves for which the curvature may vanish 
on a subset of~$I$. 
For such curves, the usual Frenet frame may not exist.
In any case, however, there exists a frame defined 
by parallel transport (\cf~\cite{Bishop_1975}).
More specifically, there exist smooth maps 
$e_2,\dots,e_{d}:I\to\Real^d$ and 
$\kappa_1,\dots,\kappa_{d-1}:I\to\Real$ 
such that $|e_\mu(s)|=1$ for all $s \in I$, 
$\mu \in \{2,\dots,d\}$
and 
\begin{equation}\label{frame}
  \begin{pmatrix}
    e_1 \\ e_2 \\ \vdots \\ e_d
  \end{pmatrix}^{\!\!\mbox{\LARGE$.$}} 
  =
  \begin{pmatrix}
    0 & \kappa_1 & \dots & \kappa_{d-1} \\
    -\kappa_1 & 0 & \dots & 0 \\
    \vdots & \vdots &  \ddots & \vdots \\ 
    -\kappa_{d-1} & 0 & \dots & 0 
  \end{pmatrix}
  \begin{pmatrix}
    e_1 \\ e_2 \\ \vdots \\ e_d
  \end{pmatrix}
  .
\end{equation}
One has $\kappa^2=\kappa_1^2+\dots+\kappa_{d-1}^2$.
Note that $\{e_1,\dots,e_d\}$ is an orthonormal 
vector field along~$\Gamma$
and that the vectors $e_2,\dots,e_d$ 
forms a basis of the normal bundle.

We locally parameterise~$\Omega_a$ by the Fermi coordinates 
\begin{align*}
\phi:\, I \times D_a &\to \Real^d \,, 
\\
(s,t) & \mapsto \gamma(s) + t_\mu e_\mu(s)
\,,
\end{align*}
where $D_a := \{t\in\Real^{d-1}:|t|<a\}$ 
is the $(d-1)$-dimensional ball, $t:=(t_2,\dots,t_d)$ 
and the Einstein summation convention is assumed,
with the range of Greek indices being $2,\dots,d$.
Using the formula~\eqref{frame},
one easily finds that the metric
$
  G := \nabla\phi \cdot (\nabla\phi)^\mathsf{T}
$
acquires the diagonal form $G=\diag(f^2,1,\dots,1)$
with the Jacobian
$$
  f(s,t) = 1 - \kappa_\mu(s) \, t_\mu
  \,.
$$
The condition that~$\Omega_a$ does not overlap itself
in particular requires that the Jacobian is positive.
A sufficient condition for the latter is that 
$a \|\kappa\|_\infty < 1$, 
where~$\|\cdot\|_\infty$ denotes the supremum norm.

Now we can compute the volume of the piece 
$\Omega_a^I := \phi(I \times D_a)$ of $\Omega_a$.
More generally, for every $r \in (0,a]$, we have
\begin{equation}\label{volume}
  |\Omega_r^I|
  = \int_I \int_{D_r} f(s,t) \, \der t \, \der s
  = \int_I \int_{D_r} 1 \, \der t \, \der s
  = |I| |D_r|
  \,,
\end{equation}
where the second equality follows by the fact 
that~$0$ is the centre of mass of~$D_r$,
so that $\int_{D_r} t \, \der t = 0$. 
Here~$|I|$ and~$|D_r|$ denote the length of the interval~$I$
and volume of the ball~$D_r$, respectively.

Let $|\partial\Omega_r^I|$ denote the $(d-1)$-dimensional
Hausdorff measure of the surface $\phi(I\times\partial D_a)$.
Expressing $|\partial\Omega_r^I|$ as the derivative of $|\Omega_r^I|$
(\cf~\cite[Lem.~3.13]{Gray}) 
and using the scaling $|D_r| = r^{d-1} |D_1|$,
we get
\begin{equation}\label{area}
  |\partial\Omega_r^I| 
  = \frac{\der}{\der r} |\Omega_r^I|
  = (d-1) r^{d-2} \, |I| |D_1|
  = \frac{d-1}{r} \, |I| |D_r|
  \,.
\end{equation}

Recalling that the curve~$\Gamma$ is parameterised 
by its arc-length via~$\gamma$, 
the local formulae~\eqref{volume} and~\eqref{area}
yield the global identities
$$
  |\Omega_a^I| = |\Gamma| |D_a|
  \qquad \mbox{and} \qquad
  |\partial\Omega_a| = \frac{d-1}{a} \, |\Gamma| |D_a|
  \,.
$$
Choosing $S:=\Omega_a$ in~\eqref{Cheeger},
we therefore get the desired upper bound
\begin{equation}\label{upper}
  h(\Omega_a) \leq  \frac{d-1}{a} 
  \,.
\end{equation}

\subsection{Proof of Theorem \ref{Thm.main}: the lower bound}\label{Sec.lower}
Now we employ a different local parameterisation of~$\Omega_a$, namely,
\begin{equation*}
\begin{aligned}
  \tilde{\phi}: I \times U &\to \Real^d \,, 
  \\
  (s,\theta) &\mapsto \gamma(s) + a \, \sigma_k(\theta) \, e_k(s)
  \,,
\end{aligned}
\end{equation*}
where $\sigma:U\to\Sphere_+^{d-1}\subset\Real^{d}$ 
is a parameterisation of the half-sphere
$
  \Sphere_+^{d-1} 
  := \{x\in\Real^d: |x|=1 \ \land \ x_1 > 0\}
$ 
and for $\theta:=(\theta_2,\dots,\theta_d) \in U$
one can choose for instance the hyperspherical coordinates.
Again, we assume the Einstein summation convention,
with the range of Latin indices being $1,\dots,d$.
Using the formula~\eqref{frame},
one easily finds that the Jacobi matrix reads
$$
  J := \nabla\tilde{\phi} =
  \begin{pmatrix}
    (1-a\,\sigma_\mu\kappa_\mu) e_1 + a \, \sigma_1 \, \kappa_\mu e_\mu \\
    a \, \partial_2\sigma_k \, e_k \\
    \vdots \\
    a \, \partial_d\sigma_k \, e_k
  \end{pmatrix} 
  ,
$$
where~$e_k$'s are arranged as row vectors. 

In order to get the same lower bound 
as the right-hand side of~\eqref{upper},
we employ Proposition~\ref{Prop.Grieser}.
Our choice of the test vector field reads, locally,
\begin{equation}\label{field}
  V(x) := \frac{x - \gamma(s)}{a} = \sigma_k(\theta) \, e_k(s)
  \,,
\end{equation}
where the relationship between $x\in\Real^d$
and $(s,\theta) \in I \times U$ is given by $x = \tilde\phi(s,\theta)$.
Clearly, $|V|=1$. It remains to compute the divergence of~$V$.

Employing the first equality of~\eqref{field}, we have
$$
  \divergence V = \frac{1}{a} \left(d-e_1\cdot\nabla s\right)
  \,,
$$ 
where~$s$ is understood as the first component of the inverse $\tilde\phi^{-1}(x)$ and the gradient acts with respect to~$x$.
Using that $|\sigma|^2=1$, 
so that $\sigma\cdot\partial_\mu\sigma=0$ for every $\mu\in\{2,\dots,d\}$,
it is straightforward to check that 
$$
  J^{-1} = 
  \begin{pmatrix}
    \displaystyle
    \frac{\sigma_k \, e_k}{\sigma_1}, *_2, \dots, *_d
  \end{pmatrix}
  ,
$$
where~$e_k$'s and~$*_\mu$'s are arranged as column vectors
and the explicit values of asterisks are not important for our purposes. 
Indeed, this formula is enough to conclude that 
$$
  \nabla s = \frac{\sigma_k \, e_k}{\sigma_1}
$$
and thus
$$
  \divergence V = \frac{d-1}{a} 
  \,.
$$ 
In view of Proposition~\ref{Prop.Grieser},
we therefore get the desired lower bound
\begin{equation}\label{lower}
  h(\Omega_a) \geq  \frac{d-1}{a} 
  \,.
\end{equation}

Recalling also~\eqref{upper},
we have just established~\eqref{Cheeger.tube}.
Recalling in addition the way how~\eqref{upper} was proved,
we also have $\mathcal{C}_{\Omega_a}=\Omega_a$. 
This concludes the proof of Theorem~\ref{Thm.main}.

\begin{Remark}[Unbounded tubes]
Let us consider the same definition~\eqref{tube}
with~$\Gamma$ being an unbounded \emph{complete} curve.  
Then the same test vector field~$V$ as in Section~\ref{Sec.lower} 
yields the lower bound~\eqref{lower}. 
At the same time, 
choosing the test domain $S:=\phi(I \times D_a)$ in~\eqref{Cheeger}
and following the (local) procedure of Section~\ref{Sec.upper},
we arrive 
\begin{equation*} 
  h(\Omega_a) \leq  
  \frac{\displaystyle\frac{d-1}{a} \, |I| |D_a| + 2 |D_a|}{|I||D_a|}
  \,.
\end{equation*}
Sending the length~$|I|$ to infinity, 
we finally get the validity of~\eqref{Cheeger.tube}
even in this unbounded case.
Note, however, that the infimum in~\eqref{Cheeger} is not achieved
and there is no Cheeger set~$\mathcal{C}_{\Omega_a}$ now.  
\end{Remark}
\subsection{Proof of Theorem \ref{Thm.shell}}
Let us denote by $\omega_d$ the volume of an Euclidean ball of radius $1$ in $\Real^d$. Since 
\[
|\partial A_{r,R}| = d\, \omega_d \, (R^{d-1}+r^{d-1})
\qquad\text{and}\qquad 
|A_{r,R}| = \omega_d \, (R^d-r^d)\,,
\]
we immediately obtain the upper bound
\[
h(A_{r,R}) 
\le \frac{|\partial A_{r,R}|}{|A_{r,R}|} 
= d \, \frac{R^{d-1}+r^{d-1}}{R^d-r^d}\,.
\]

In order to conclude the proof we only need to show the opposite inequality and discuss the equality case. As before, this will be accomplished by applying Proposition \ref{Prop.Grieser} to a suitable vector field $V$, that is now defined as
\[
V(x) := f(|x|)\,x\,,
\]
where
$$
  f(t) := C -(Cr^d+r^{d-1}) \, t^{-d}
  \qquad \mbox{with} \qquad
  C:=\frac{R^{d-1}+r^{d-1}}{R^d-r^d}
  \,.
$$ 
Notice that the function $tf(t)$ is strictly increasing when $t>0$, hence by easy computations we obtain $Rf(R) = -rf(r) = 1$ and
\[
t\,|f(t)| < 1 \,, \qquad \forall\, t\in (r,R)\,.
\]
Then $|V(x)| <1$ for all $r<|x|<R$, and moreover $V(x)$ agrees with the exterior unit normal on $\partial A_{r,R}$. On the other hand, we have
\[
\divergence V(x) 
= f'(|x|)|x| + d\,f(|x|) 
= d\,(Cr^d+r^{d-1})|x|^{-d} + d\,\big(C-(Cr^d+r^{d-1})|x|^{-d}\big) 
= dC 
= \frac{|\partial A_{r,R}|}{|A_{r,R}|}\,,
\]
hence by Proposition \ref{Prop.Grieser} we obtain the required lower bound 
\[
h(A_{r,R}) \ge d \, \frac{R^{d-1}+r^{d-1}}{R^d-r^d}\,.
\]

The uniqueness of the Cheeger set of $A_{r,R}$ 
follows from the fact that $|V(x)|$ 
is strictly less than~$1$ in $A_{r,R}$, 
hence the strict inequality
\[
|\partial A| > h(A_{r,R}) \, |A|
\]
holds for any measurable set $A\subset A_{r,R}$, 
for which both $|A|$ and $|A_{r,R}\setminus A|$ are strictly positive. 
This concludes the proof of Theorem \ref{Thm.shell}.

\subsection*{Acknowledgment}
%
The research of D.K.\ was partially supported 
by the GACR grant No.\ 18-08835S
and by FCT (Portugal) through project PTDC/MAT-CAL/\-4334/\-2014.
 

%
%

{\footnotesize
\providecommand{\bysame}{\leavevmode\hbox to3em{\hrulefill}\thinspace}
\providecommand{\MR}{\relax\ifhmode\unskip\space\fi MR }
\providecommand{\MRhref}[2]{%
  \href{http://www.ams.org/mathscinet-getitem?mr=#1}{#2}
}
\providecommand{\href}[2]{#2}

}
\end{document}